\documentclass[a4paper,12pt]{amsart}
\usepackage[letterpaper, margin=1.2in]{geometry}
\usepackage{latexsym}
\usepackage{amssymb}
\usepackage{amsmath}
\usepackage{amsfonts}
\usepackage[table]{xcolor}
\usepackage{pgfplots}
\usepackage{color}
\usepackage{multirow}
\usepackage{graphicx}
\usepackage{txfonts,pxfonts} 
\usepackage{tikz}
\usetikzlibrary{calc,arrows}
\usetikzlibrary{calc}
\usepackage{verbatim}
\everymath{\displaystyle}
\newtheorem{thm}{Theorem}[section]

\newtheorem{lem}[thm]{Lemma}

\newtheorem{cor}[thm]{Corollary}

\newtheorem{prop}[thm]{Proposition}

\theoremstyle{definition}

\newtheorem*{examples}{Examples}

\newcommand{\Z}{\mathbb Z}

\newcommand{\Q}{\mathbb Q}

\def\ol{\overline}
\def\al{\alpha}
\def\la{\lambda}

\def\th{\theta}
\def\md#1{\ \mbox{\rm(mod }{#1})}

\newcounter{cs}
\stepcounter{cs}
\newcommand{\casos}{\begin{itemize}}
	\newcommand{\fcasos}{\end{itemize}\setcounter{cs}{1}}

\newfont{\tit}{cmr12 scaled \magstep3}
\title{On  index  divisors and monogenity of certain sextic number fields defined by $x^6+ax^5+b$}
\author{Lhoussain El Fadil}
\author{Omar Kchit}
\address{Faculty of Sciences Dhar El Mahraz, P.O. Box  1796 Atlas-Fes, Sidi Mohamed ben Abdellah University,  Morocco}
\email{lhoussain.elfadil@usmba.ac.ma, \,\, orcid: 0000.0003.4175.8064} 
\email{omar.kchit@usmba.ac.ma, \,\, orcid: 0000.0002.0844.5034}
\begin{document}
		\maketitle

		{
		\begin{abstract}
     The main goal of  this paper is to provide  a complete answer to the  Problem 22 of  Narkiewicz \cite{Nar} for any sextic number  field $K$ generated  by a complex  root  $\al$ of a monic irreducible trinomial  $F(x) = x^6+ax^5+b \in \Z[x]$. Namely we   calculate the index of the field $K$. In particular, if $i(K)\neq 1$, then $K$ is not mongenic. Finally, we illustrate our results by some computational examples.        
\end{abstract}
\section{Introduction} 
Let $K$ be an number field and $\Z_K$  its ring of integers. For any primitive element $\al\in \Z_K$ of $K$, let
$ind(\al) = (\Z_K: \Z[\al])$ be the index of $\Z[\al]$ in $\Z_K$. The greatest common divisor of the indices of all
integral primitive elements of $K$ is called the  index of $K$, and  denoted by $i(K)$. Say $i(K)=\gcd \ \{ ( \Z_K; \Z [\alpha]) \, |\,\alpha \in \widehat{\Z}_K  \}$. A prime $p$ dividing $i(K)$ is called a prime common index divisor of $K$. If $\mathbb{Z}_K$ has a power integral basis, then the index of $K$ is trivial, namely $i(K)=1$. Therefore a field having a prime common index divisor is not monogenic. 
Problem 22 of Narkiewicz \cite{Nar} asks for an explicit formula for the highest power of a
given prime $p$ dividing $i(K)$, say $\nu_p(i(K))$.
 The first  number field with non trivial index was given by Dedekind in 1871,
 	 who exhibited examples in cubic and quartic number fields. For example, he considered the cubic field $K$ generated by a complex root of $x^3-x^2-2x-8$ and  showed that the prime $2$ splits completely in $K$. So, if we suppose that $K$ is monogenic, then we would be able to find a cubic polynomial generating $K$, that splits completely into distinct polynomials of degree $1$ in $\mathbb{F}_2[x]$. Since there are only $2$ distinct polynomials of degree $1$ in $\mathbb{F}_2[x]$, this is impossible. Based on these ideas and using Kronecker's theory of algebraic number fields, Hensel  gave necessary and sufficient conditions on the so-called "index divisors of $K$" for any prime $p$ to be  a prime common index divisor \cite{He2}.
In \cite{En}, for any number field of degree $n\le 7$, Engstrom gave an explicit formula which relates   $\nu_p(i(K))$ to the factorization of  $p\mathbb{Z}_K$ into powers of prime ideals of $K$ for every positive prime $p\le n$.
In \cite{Nak}, Nakahara  studied the  index of non-cyclic but abelian biquadratic number fields. 
 In \cite{GPP}  Ga\'al et al. characterized  the field indices of biquadratic number fields having Galois group $V_4$. In \cite{DS} for any  quartic number field $K$ defined by a trinomial $x^4 + ax + b$, Davis and  Spearman gave necessary and sufficient conditions on $a$ and $b$, which characterize when $2$ divides $i(K)$.  
Also in   \cite{EG} for any  quartic number field $K$ defined by a trinomial $x^4 + ax^2 + b$, El Fadil and Ga\'al gave necessary and sufficient conditions on $a$ and $b$, which characterize when a prime $p$ divides $i(K)$.
In \cite{E5} for any  quintic number field $K$ defined by a trinomial $x^5 + ax^2 + b$, we gave necessary and sufficient conditions on $a$ and $b$, which characterize when a prime  $p$ divides $i(K)$.
In this paper, for any sextic number field $K$ defined by a trinomial $x^6 + ax^5 + b\in \Z[x]$, we characterize the index of $K$. Based on Engstrom's results given in \cite{En},  the unique prime candidates to divide $i(K)$ are $2$, $3$ and $5$; $i(K)=2^{\nu_2}3^{\nu_3}5^{\nu_5}$ , with $0\le \nu_2\le 5$, $0\le \nu_3\le 2$, and $0\le \nu_5\le 2$.  In Theorems \ref{thm2} and \ref{thm3} and Proposition \ref{pro5}, we evaluate $\nu_p=\nu_p(i(K))$ for every prime $p=2,3,\text{ and }5$. We show that $0\le \nu_2\le 2$, $0\le \nu_3\le 1$ and $\nu_5=0$. Namely $i(K)\in\{1,2,3,4,6,12\}$.}
\section{main results}
{Let $a$ and $b$ be two rational integers such that $F(x)=x^6+ax^5+b$ is irreducible over $\Q$, and let $K = \Q(\al)$ be a number field generated by a root $\al$ of  $F(x)$.
 If $\nu_p(a)\ge 1$ and $\nu_p(b)\ge 6$ for a prime $p$, then $\th=\frac{\al}{p}$ is a primitive integral element of $K$, which is root of the irreducible 
trinomial $x^6 + Ax^5 +B$  having integer coefficients $A=a/p$ and $B=b/p^6$. So up to replace $\al$ by $\theta$ and  repeat this process until to get either $\nu_p(a)< 1$ or $\nu_p(b)< 6$, we can assume that $\nu_p(a)< 1$ or $\nu_p(b)< 6$ for every prime $p$. Let $\Z_K$ be the ring of integers of $K$, $\triangle(F)=-b^4(6^6\times b-5^5\times a^6)$ the discriminant of $F(x)$ and $d_K$ the absolute discriminant of $K$. Then we have the following well know formula: 
 \begin{eqnarray}\label{indexdisc}
 \triangle(F)=\pm ind(\al)^2d_K,
\end{eqnarray} which relates the index $ind(\al)$, $\triangle(F)$,  and $d_K$.}

	The following theorem is a refinement of \cite[Theorem 3.4]{Smt}.
	\begin{thm}\label{thm1}
		The ring $\Z[\alpha]$ is the ring of integers of $K$ if and only if { the following conditions hold}:
		\begin{enumerate}
			\item[(i)] {$b$ is square free}.
			\item[(ii)] {If $2\mid a\text{ and }2\nmid b,\text{ then }(a,b)\in\left\{(0,1),(2,3)\right\}\md{4}$}.
			\item[(iii)] {If $3\mid a \text{ and }3\nmid b, \text{ then }(a,b)$ is in the set} $$\left\{(0,2),(0,4),(0,5),(0,7),(3,-1),(3,1),(3,4),(3,7),(6,-1),(6,1),(6,4),(6,7)\right\}\md{9}.$$
			\item[(iv) ]{For every {prime} $p$ dividing $(6^6b-5^5a^6)$, if $p\nmid ab,\text{ then }p^2\nmid (6^6b-5^5a^6)$}.
		\end{enumerate}
In particular, if these conditions hold, then  the index of $K$ is trivial, say $i(K)=1$.
	\end{thm}

	{In the remainder of this section, for every prime $p$, we give sufficient and necessary conditions on $a$ and $b$ so that  $p$ divides $i(K)$ for every prime $p=2,3$. Furthermore in every case, $\nu_p(i(K))=\nu_p$ will be calculated for $p=2,3$.}
\smallskip

	Recall that, for a {prime} $p$ and a monic polynomial $\phi\in\mathbb{Z}[x]$ whose reduction is irreducible over $\mathbb{F}_p$, let $N_{\phi}^+(F)=S_1+\cdots+S_t$ be the principal $\phi$-Newton polygon of $F(x)$ with respect to $p$. We say that $F(x)$ is a $\phi$-regular polynomial, if the residual polynomial $R_{\lambda}(F)(y)$ attached to each side of $N_{\phi}^+(F)$ is separable over $\mathbb{F}_{\phi}$.\\
	\smallskip
	\newline
	
	Let $p$ be a {prime} and $K$ a number field defined by $F(x)\in\mathbb{Z}[x]$. Let $\phi=x-z\in\mathbb{Z}[x]$ and {$F(x)=\phi^6+a_5\phi^5+a_4\phi^4+a_3\phi^3+a_2\phi^2+a_1\phi+a_0$, for some $(a_0,\dots,a_5)\in\mathbb{Z}^6$}. Assume that $N_{\phi}^+(F)$ has a side of integer slope $-k$ $(k\in\mathbb{N})$ {and length} $l\geq 2$. Assume also that the residual polynomial {$R_{\lambda}(F)(y)$} attached to this side has a double root $\ol{t}\in \mathbb{F}_{\phi}\simeq\mathbb{F}_p$. Then we can construct an element $s\in\mathbb{Z}$, such that $s\equiv z\md{p}$ and $F(x)$ is $(x-s)$-regular with respect to $p$. Such an element $s$ is called a regular integer of $F(x)$ with respect to $\ol{\phi}$ and $p$. Indeed, by theorems of the polygon and of the residual polynomial, $F(x)=F_1(x)F_2(x)$ in $\mathbb{Z}_p[x]$, such that $F_2$ is monic, $N_{\phi}(F_2)$ has a single side of slope $-k$, and $R_{\lambda}(F_2)(y)=c(y-t)^2$, {with $c\in\mathbb{F}_{\phi}^*$}, is the residual polynomial attached to this side. Let {$s_0=z$}, $s_1=s_0+p^kt$. Then $F_2(x)=b_2(x-s_1)^2+b_1(x-s_1)+b_0$ for some $(b_0,b_1,b_2)\in\mathbb{Z}_p^3$, such that $\nu_p(b_2)=0$, $\nu_p(b_1)\geq k+1$, and $\nu_p(b_0)\geq 2k+1$. If $N_{x-s_1}(F_2)$ has a slope $-h$ for some integer $h\geq k$ and {$R_{\lambda}(F_2)(y)=c_1(y-t_1)^2$} for some integer $t_1\in \mathbb{Z}$ {and $c_1\in\mathbb{F}_{\phi}^*$}, then we can repeat the same process. In this case, we have $ind_{x-s_1}(F)\geq ind_{\phi}(F)+1$. Thus $(\mathbb{Z}_K:\mathbb{Z}[\alpha])\geq ind_{x-s_1}(F)\geq ind_{\phi}(F)+1$. Since $(\mathbb{Z}_K:\mathbb{Z}[\alpha])$ is finite, this process cannot continue infinitely. Thus after a finite
	number of iterations, this process will provide a regular element with respect to {$\overline{\phi}$} and $p$.\\
	\smallskip
	
	\begin{thm}\label{thm2} {Let $u=\nu_2(5a+6)$ and $v=\nu_2(a+b+1)$}.
		The {prime} $2$ {divides the index $i(K)$} if and only if one of the following conditions holds:
		\begin{enumerate}
			\item {$(a,b)\in\{(0,3),(0,-1)\}\md{8}$}. {In this case, $\nu_2(i(K))=2$}.
			\item {$a\equiv 2\md{4}$ and $b\equiv -(a+1)+2^{2u}\md{2^{2u+1}}$}. {In this case, $\nu_2(i(K))=1$}.
			\item {$(a,b)\equiv (2,1)\md{4}$, $v<2u$, $v=2k$ for some positive integer $k$}, and
			$\nu_2(b+as^5+s^6)=2\nu_2(5as^4+6s^5)$, for some integer $s$ such that $s\equiv 1\md{2}$ and $F(x)$ is $(x-s)$-regular with respect to $2$. {In this case, $\nu_2(i(K))=1$}.
		\end{enumerate}
		In particular, if one of these conditions holds, then $K$ is not monogenic.
	\end{thm}
	\smallskip
	Before to state Theorem $\ref{thm3}$,
	for every integer $m\in\mathbb{Z}$, {let} $m_3=\frac{m}{3^{\nu_3(m)}}$.
	\begin{thm}\label{thm3}
		Let ${u}=\nu_2(5a+6)$, ${v}=\nu_3(a+b+1)$, ${\mu}=\nu_3(5a-6)$, and ${\tau}=\nu_3(-a+b+1)$. 
		The {prime $3$ divides the index $i(K)$} if and only if one of the following conditions holds:
		\begin{enumerate}
			\item {$(a,b)\equiv (0,-1)\md{9}$}.
			\item $a\equiv 12\md{27}$, $b\equiv a-1\md{81}$, and
			$6\leq2{\mu}<{\tau}+1$ or $5\leq {\tau}+1<2{\mu},~{\tau}$ is odd, and $(-a+b+1)_3\equiv-1\md{3}$.
			
			\item $(a,b)\in\{(21,74),(48,20),(75,47)\}\md{81}$ and 
			$6\leq 2\nu_3(5as^4+6s^5)<\nu_3(b+as^5+s^6)+1$ or $5\leq\nu_3(b+as^5+s^6)+1<2\nu_3(5as^4+6s^5),~\nu_3(b+as^5+s^6)$ is odd, and $(b+as^5+s^6)_3\equiv-1\md{3}$, where $s$ is an integer such that $F(x)$ is $(x-s)$-regular with respect to $3$ and $s\equiv -1\md{3}$.
			
			\item $a\equiv 3\md{27}$, $b\equiv a-1\md{81}$, and
			$6\leq 2\nu_3(5as^4+6s^5)<\nu_3(b+as^5+s^6)+1$ or $5\leq\nu_3(b+as^5+s^6)+1<2\nu_3(5as^4+6s^5),~\nu_3(b+as^5+s^6)$ is odd, and $(b+as^5+s^6)_3\equiv-1\md{3}$, where $s$ is an integer such that $F(x)$ is $(x-s)$-regular with respect to $3$ and $s\equiv -1\md{3}$.
			
			\item $a\equiv 12\md{27}$, $b\equiv a-1\md{81}$, $2{\mu=\tau}+1$, $(-a+b+1)_3\equiv 1\md{3}$, {$(5a-6)_3\equiv \pm1\md{3}$}, and $6\leq 2\nu_3(5as^4+6s^5)<\nu_3(b+as^5+s^6)+1$ or $5\leq\nu_3(b+as^5+s^6)+1<2\nu_3(5as^4+6s^5),~\nu_3(b+as^5+s^6)$ is odd, and $(b+as^5+s^6)_3\equiv-1\md{3}$, where $s$ is an integer such that $F(x)$ is $(x-s)$-regular with respect to $3$ and $s\equiv -1\md{3}$.

			\item $a\equiv 15\md{27}$, $b\equiv -a-1\md{81}$ and $6\leq2{u}<{v}+1$ or $5\leq{v}+1<2{u},~{v}$ is odd, and $(a+b+1)_3\equiv -1\md{3}$.
			
			\item $(a,b)\in\{(6,47),(33,20),(60,74)\}\md{81}$ and 
			$6\leq 2\nu_3(5as^4+6s^5)<\nu_3(b+as^5+s^6)+1$ or $5\leq\nu_3(b+as^5+s^6)+1<2\nu_3(5as^4+6s^5),~\nu_3(b+as^5+s^6)$ is odd, and $(b+as^5+s^6)_3\equiv-1\md{3}$, where $s$ is an integer such that $F(x)$ is $(x-s)$-regular with respect to $3$ and $s\equiv 1\md{3}$.
			
			\item $a\equiv24\md{27}$, $b\equiv a-1\md{81}$, and $6\leq 2\nu_3(5as^4+6s^5)<\nu_3(b+as^5+s^6)+1$ or $5\leq\nu_3(b+as^5+s^6)+1<2\nu_3(5as^4+6s^5),~\nu_3(b+as^5+s^6)$ is odd, and $(b+as^5+s^6)_3\equiv-1\md{3}$, where $s$ is an integer such that $F(x)$ is $(x-s)$-regular with respect to $3$ and $s\equiv 1\md{3}$.
			
			\item $a\equiv 15\md{27}$, $b\equiv -a-1\md{81}$, $2{u=v}+1$, $(a+b+1)_3\equiv1\md{3}$, {$(5a+6)_3\equiv \pm1\md{3}$}, and $6\leq 2\nu_3(5as^4+6s^5)<\nu_3(b+as^5+s^6)+1$ or $5\leq\nu_3(b+as^5+s^6)+1<2\nu_3(5as^4+6s^5),~\nu_3(b+as^5+s^6)$ is odd, and $(b+as^5+s^6)_3\equiv-1\md{3}$, where $s$ is an integer such that $F(x)$ is $(x-s)$-regular with respect to $3$ and $s\equiv 1\md{3}$.
		\end{enumerate}
		
		In particular, if one of the above conditions holds, then $K$ is not monogenic.\\
		{Furthermore, $\nu_3(i(K))=1$ in all the above cases}.
	\end{thm}
	\begin{prop}\label{pro5}
		For every $(a,b)\in\mathbb{Z}^2$, {$5$ does not divide the index $i(K)$}.
	\end{prop}
	\smallskip
\begin{cor}\label{cor}
		Let $K$ be a sextic number field generated by a monic irreducible polynomial $F(x)=x^6+ax^5+b\in \mathbb{Z}[x]$, then the following statements hold:
		\begin{enumerate}
			\item If $a\equiv 0\md{72}$ and $b\notin\{ -1,3,7,8,11,15,17,19,23,26,27,31,35,39,42,43,44,47$,\\
			$51,53,55,59,62,63,67\}\md{72}$ , then by Theorem \ref{thm2} (1), $\nu_2(i(K))=0$ and by Theorem \ref{thm3} (1), $\nu_3(i(K))=0$. Thus $i(K)=1$.
			
			\item If $(a,b)\in\{(14,1),(14,33),(22,25),(22,57),(38,9),(38,41),(46,1),(46,33),(62,17)$,\\
			$(62,49)\}\md{96}$, then by Theorem \ref{thm2} (2), $\nu_2(i(K))=1$ and by Theorem \ref{thm3}, $\nu_3(i(K))=0$. Thus $i(K)=2$.
			
			\item If $a\equiv 0\md{72}$ and $b\equiv 8,17,26,44,62\md{72}$, then by Theorem \ref{thm2} (1), $\nu_2(i(K))=0$ and by Theorem \ref{thm3} (1), $\nu_3(i(K))=1$. Thus $i(K)=3$.
			
			\item If $a\equiv 0\md{72}$ and $b\equiv 3,7,11,15,19,23,27,31,39,43,51,55,59,63,67\md{72}$, then by Theorem \ref{thm2} (1), $\nu_2(i(K))=2$ and by Theorem \ref{thm3} (1), $\nu_3(i(K))=0$. Thus $i(K)=4$.
			
			\item If $(a,b)\in\{(54,89),(126,17),(198,233),(270,161)\}\md{288}$, then by Theorem \ref{thm2} (2), $\nu_2(i(K))=1$ and by Theorem \ref{thm3} (1), $\nu_3(i(K))=1$. Thus $i(K)=6$.
			
			\item If $(a,b)\in\{(0,-1),(0,35)\}\md{72}$, then by Theorem \ref{thm2} (1), $\nu_2(i(K))=2$ and by Theorem \ref{thm3} (1), $\nu_3(i(K))=1$. Thus $i(K)=12$.
			
			
			
			
			
		\end{enumerate}
	\end{cor}
	
	\color{black}
	\section{Preliminaries}
	Newton polygon techniques is a standard method which is rather technical but very efficient to apply. We briefly describe the use of these techniques, which makes our proofs understandable. For more details, we refer to \cite{EMN} and \cite{GMN}.\\
	Let $K=\mathbb{Q}(\alpha)$ be a number field generated by a complex root $\alpha$ of a monic irreducible polynomial $F(x)\in\mathbb{Z}[x]$. We shall use Dedekind’s theorem \cite[25, Chapter I, Proposition 8.3]{Neu} relating the prime
	ideal factorization of $p\mathbb{Z}_K$ and the factorization of $F(x)$ modulo $p$ $($for primes $p$ not
	dividing $(\mathbb{Z}_K:\mathbb{Z}[\alpha]))$. Also, we shall need Dedekind’s criterion \cite[Theorem 6.1.4]{Co}
	on the divisibility of $(\mathbb{Z}_K:\mathbb{Z}[\alpha])$ by primes $p$.\\
	For any prime $p$, let $\nu_p$ be the $p$-adic valuation of $\mathbb{Q}$, $\mathbb{Q}_p$ its $p$-adic completion, and $\mathbb{Z}_p$ the ring of $p$-adic integers. We denote by $\nu_p$ be the Gauss's extension of $\nu_p$ to $\mathbb{Q}_p(x)$ and is defined on $\mathbb{Q}_p[x]$ by $\nu_p(\sum_{i=0}^na_ix^i)=\min_i\{\nu_p(a_i)\},a_i\in\mathbb{Q}_p$. Also, for nonzero polynomials, $P,Q\in\mathbb{Q}_p[x]$, we extend this valuation to $\nu_p(P/Q)=\nu_p(P)-\nu_p(Q)$. Let $\phi\in\mathbb{Z}_p[x]$ be a monic lift to an irreducible factor of $F(x)$ modulo $p$. Upon the Euclidean division by successive powers of $\phi$, there is a unique {$\phi$-expansion} of $F(x)$; that is $F(x)=a_0(x)+a_1(x)\phi(x)+\cdots+a_l(x)\phi(x)^l$, where $a_i(x)\in\mathbb{Z}_p[x]$ and deg$(a_i)<$deg$(\phi)$. For every $i=0,\dots,l$, let  $u_i=\nu_p(a_i(x))$. The {$\phi$-Newton polygon} of $F(x)$ with respect to $p$, is the lower boundary convex envelope of the set of points $\{(i,u_i),~ a_i(x)\neq 0\}$ in the Euclidean plane, which we denote by $N_{\phi}(F)$. It is the process of joining the obtained edges $S_1,\dots,S_t$ ordered by increasing slopes, which can be expressed as $N_{\phi}(F)= S_1+\cdots+S_t$. For every
	side $S_i$ of $N_{\phi}(F)$, its length $l(S_i)$ is the length of its projection to the $x$-axis and its
	height $h(S_i)$ is the length of its projection to the $y$-axis. We call $d(S_i)=$gcd$(l(S_i),h(S_i))$ the degree of $S_i$. The polygon determined by the sides of the $\phi$-Newton polygon with negative slopes is called the {principal $\phi$-Newton polygon} of $F(x)$, and it is denoted by $N_{\phi}^+(F)$. As defined in \cite[Def. 1.3]{EMN}, the $\phi$-{index} of $F(x)$, denoted  $ind_{\phi}(F)$, is  deg$(\phi)$ multiplied by the number of points with natural integer coordinates that lie below or on the polygon $N_{\phi}^+(F)$, strictly above the horizontal axis,{ and strictly beyond the vertical axis}. Let $\mathbb{F}_{\phi}$ be the field $\mathbb{F}_p[x]/(\ol{\phi})$, then to every side $S$ of $N_{\phi}^+(F)$, with initial point $(i,u_i)$, and every $i=0,\ldots,l$, let the residue coefficient $c_i\in\mathbb{F}_{\phi}$ defined as follows: 
	$$c_{i}=
	\left
	\{\begin{array}{ll} 0,& \mbox{ if } (s+i,{\it u_{s+i}}) \mbox{ lies strictly
			above } S,\\
		\left(\dfrac{a_{s+i}(x)}{p^{{\it u_{s+i}}}}\right)
		\,\,
		\mod{(p,\phi(x))},&\mbox{ if }(s+i,{\it u_{s+i}}) \mbox{ lies on }S.
	\end{array}
	\right.$$
	where $(p,\phi(x))$ is the maximal ideal of $\mathbb{Z}[x]$ generated by $p$ and $\phi$. Let $\lambda=-h/e$ be the slope of $S$, where  $h$ and $e$ are two positive coprime integers and $l=l(S)$. Then  $d=l/e$ is the degree of $S$. Since  the points  with integer coordinates lying{ on} $S$ are exactly ${(s,u_s),(s+e,u_{s}-h),\dots, (s+de,u_{s}-dh)}$. Thus if $i$ is not a multiple of $e$, then 
	$(s+i, u_{s+i})$ does not lie on $S$, and so $c_i=0$. Let
	${R_{\lambda}(F)(y)}=t_dy^d+t_{d-1}y^{d-1}+\cdots+t_{1}y+t_{0}\in\mathbb{F}_{\phi}[y]$, called  
	the {residual polynomial} of $F(x)$ associated to the side $S$, where for every $i=0,\dots,d$,  $t_i=c_{ie}$. If ${R_{\lambda}(F)(y)}$ is square free for each side of the polygon $N_{\phi}^+(F)$, then we say that $F$ is $\phi$-regular. Let $\ol{F(x)}=\prod_{i=1}^{r}\ol{\phi_i}^{l_i}$ be the factorization of $F(x)$ into powers of monic irreducible coprime polynomials over $\mathbb{F}_p$, we say that the polynomial $F(x)$ is $p$-{regular} if $F(x)$ is a $\phi_i$-regular polynomial with respect to $p$ for every $i=1,\dots,r$. Let  $N_{\phi_i}^+=S_{i1}+\cdots+S_{ir_i}$ be the $\phi_i$-principal Newton polygon of $F(x)$ with respect to $p$. For every $j=1,\dots,r_i$, let let $R_{\lambda_{ij}}(y)=\prod_{s=1}^{s_{ij}}\psi_{ijs}^{a_{ijs}}(y)$ be the factorization of $R_{\lambda_{ij}}(y)$ in $\mathbb{F}_{\phi_i}[y]$. Then we have the following  theorem of index of Ore \cite[Theorem 1.7 and Theorem 1.9]{EMN} and \cite{O, MN}:
	\begin{thm}\label{thm4}
		Under the above hypothesis, we have the following:
		\begin{enumerate}
			\item 
			$$\nu_p((\mathbb{Z}_K:\mathbb{Z}[\alpha]))\geq\sum_{i=1}^{r}ind_{\phi_i}(F).$$  
			The equality holds if $F(x) \text{ is }p$-regular.
			\item 
			If  $F(x) \text{ is }p$-regular, then
			$$p\mathbb{Z}_K=\prod_{i=1}^r\prod_{j=1}^{t_i}\prod_{s=1}^{s_{ij}}\mathfrak{p}_{ijs}^{e_{ij}}.$$
			where $e_{ij}$ is the smallest positive integer satisfying $e_{ij}\la_{ij}\in \Z$ and\\
			$f_{ijs}=\deg(\phi_i)\times \deg(\psi_{ijs})$ is the residue degree of $\mathfrak{p}_{ijs}$ over $p$ for every $(i,j,s)$.
		\end{enumerate}
	\end{thm}
\begin{cor}\label{cor1}
	Under the hypothesis  above {$($Theorem $\ref{thm4})$}, if for every $i=1,\dots,r,\,l_i=1\text{ or }N_{\phi}^+(F)=S_i$ has a single side of height $1$, then $\nu_p((\mathbb{Z}_K:\mathbb{Z}[\alpha]))=0$.
\end{cor}
	\section{Proofs of Main Results}
	{Throughout this section, $\mathfrak{p}$ is a prime ideal of $\mathbb{Z}_K$ lying above the prime $p$, $f=f(\mathfrak{p})$ its residue degree, and for every integer $m$, we denote $m_p=\cfrac{m}{p^{\nu_p(m)}}$}.\\
	{\it Proof of Theorem \ref{thm1}}.\\
	For the proof of this theorem, let us use theorem of index with "if and only if" as it is given in  \cite[Theorem 4.18]{GMN}.
	Since $\Delta(F)=-b^4(6^6b-5^5a^6)$, thanks to the  formula $\nu_p(\Delta(F))=\nu_p(d_K)+2\nu_p(ind(\alpha))$, $\mathbb{Z}[\alpha]$ is the ring of integers of $K$ if and only if $p$ does not divide $(\mathbb{Z}_K:\mathbb{Z}[\alpha])$ for {prime} $p$ dividing {$b \times(6^6b-5^5a^6)$}.
	\begin{enumerate}
		\item[(i)]
		Let $p$ be a {prime} dividing $b$. If $p\mid a$, then $F(x)\equiv x^6\md{p}$. Let $\phi=x$. So, $ind_1(F)=0$ if and only if $N_{\phi}^+(F)=S$ has a single side of height $1$, which is equivalent to $\nu_p(b)=1$, that means $b$ is square free; $p$ does not divide $(\mathbb{Z}_K:\mathbb{Z}[\alpha])$.
		If $p\nmid a$, $F(x)\equiv x^5(x+a)\md{p}$. Let $\phi_1(x)=x$ and $\phi_2(x)=x+a$. The $\phi_2$-expansion of $F(x)$ is given by
		\begin{equation}\label{eq2}
			F(x)=\phi_2^6-5a\phi_2^5+10a^2\phi_2^4-10a^3\phi_2^3+5a^4\phi_2^2-a^5\phi_2+b.
		\end{equation}
		So, $ind_1(F)=ind_{\phi_1}(F)+ind_{\phi_2}(F)=0$ if and only if $N_{\phi_i}^+(F)=S_i$ have a single side of height $1$ for every $i=1,2$, which is equivalent to $\nu_p(b)=1$, that means $b$ is square free; $p$ does not divide $(\mathbb{Z}_K:\mathbb{Z}[\alpha])$. 
		\item[(ii)] If $2\mid a\text{ and }2\nmid b$, then $F(x)\equiv x^6+1\equiv(x^2+x+1)^2(x+1)^2\md{2}$. Let $\phi_1(x)=x+1$ and $\phi_2(x)=x^2+x+1$. The $\phi_1$-expansion of $F(x)$ is given by 
		\begin{equation}\label{eq3}
			F(x)=\phi_1^6+(a-6)\phi_1^5+(-5a+15)\phi_1^4+(10a-20)\phi_1^3+(-10a+15)\phi_1^2+(5a-6)\phi_1-a+b+1.
		\end{equation}
		Thus, $ind_{\phi_1}(F)=0$ if and only if $\nu_2(-a+b+1)=1$. The $\phi_2$-expansion of $F(x)$ is given by 
		\begin{equation}\label{eq4}
			F(x)=\phi_2^3+((a-3)x-2a)\phi_2^2+((a+2)x+3a-2)\phi_2-ax-a+b+1.
		\end{equation}
		Thus, $ind_{\phi_2}(F)=0$ if and only if $\nu_2(-a+b+1)=1$ or $\nu_2(a)=1$.\\
		We conclude that $ind_1(F)=ind_{\phi_1}(F)+ind_{\phi_2}(F)=0$ if and only if $\nu_2(-a+b+1)=1$; which means that $\nu_2((\mathbb{Z}_K:\mathbb{Z}[\alpha]))=0$ if and only if $(a,b)\in\{(0,1),(2,3)\}\md{4}$.
		\item[(iii)] If $3\mid a\text{ and }3\nmid b$, then two cases arise:\\
		If $b\equiv1\md{3}$, then $F(x)\equiv x^6+1\equiv (x^2+1)^3\md{3}$. Let $\phi(x)=x^2+1$. Let 
		\begin{equation}\label{eq5}
			F(x)=\phi^3+(ax-3)\phi^2+(-2ax+3)\phi+ax+b-1.
		\end{equation}
		Thus, $ind_1(F)=0$ if and only if $\nu_3(ax+b-1)=1$; $b\equiv 4,7\md{9}\text{ or }a\equiv3,6\md{9}$, which means  $$(a,b)\in\{(0,4),(0,7),(3,1),(3,4),(3,7),(6,1),(6,4),(6,7)\}\md{9}.$$
		If $b\equiv2\md{3}$, then $F(x)\equiv (x+1)^3(x-1)^3\md{3}$. Let $\phi_1(x)=x+1$ and $\phi_2(x)=x-1$. The $\phi_1$-expansion of $F(x)$ is given above by $(\ref{eq3})$ and $ind_{\phi_1}(F)=0$ if and only if $\nu_3(-a+b+1)=1$. The $\phi_2$-expansion of $F(x)$ is given by
		\begin{equation}\label{eq6}
			F(x)=\phi_2^6+(a+6)\phi_2^5+(5a+15)\phi_2^4+(10a+20)\phi_2^3+(10a+15)\phi_2^2+(5a+6)\phi_2+a+b+1.
		\end{equation}
		Thus, $ind_{\phi_2}(F)=0$ if and only if $\nu_3(a+b+1)=1$. We conclude that $ind_1(F)=ind_{\phi_1}(F)+ind_{\phi_2}(F)=0$ if and only if $\nu_3(-a+b+1)=\nu_3(a+b+1)=1$, which means that $(a,b)\in\{(0,2),(0,5),(3,-1),(6,-1)\}\md{9}$.\\
		We conclude that if $3\mid a\text{ and }3\nmid b$, {$\nu_3((\mathbb{Z}_K:\mathbb{Z}[\alpha]))=0$} if and only if $(a,b)$ is in the set
		$$\{(0,2),(0,4),(0,5),(0,7),(3,-1),(3,1),(3,4),(3,7),(6,-1),(6,1),(6,4),(6,7)\}\md{9}.$$
		\item[(iv)] {Let $p$ be a {prime} dividing $(6^6b-5^5a^6)\text{ such that }p\nmid ab$. Thanks to the index formula above, if $p^2\nmid(6^6b-5^5a^6)$, then $p\nmid(\mathbb{Z}_K:\mathbb{Z}[\alpha])$. Conversely, if $p\nmid (\mathbb{Z}_K:\mathbb{Z}[\alpha])$, then by \cite[Theorem $1$]{LNV}, we have:
			\begin{equation}\tag{1.1}
				\nu_p(d_K)=\left\{
				\begin{array}{ll}
					0&\text{ if }\nu_p(6^6b-5^5a^6)\text{ is even},\\
					1&\text{ if }\nu_p(6^6b-5^5a^6)\text{ is odd}.
				\end{array}
				\right.	
			\end{equation}
			Then, according to the index formula and $(1.1)$, $p^2\nmid (6^6b-5^5a^6)$}.
	\end{enumerate}
	\begin{flushright}
		$\square$
	\end{flushright}
	\smallskip
	{For the proof of Theorems \ref{thm2} and \ref{thm3} and Proposition \ref{pro5}, we need the following lemma, which characterizes the prime divisors of the index $i(K)$}.
		\begin{lem}\label{index}$($\cite[Theorem 2.2]{HS}$)$\\
		Let $p$ be a prime and $K$ be a number field. For every positive integer $f$, let $\mathcal{P}_f$ be the number of distinct prime ideals of $\mathbb{Z}_K$ lying above $p$ with residue degree $f$ and $\mathcal{N}_f$ the number of monic irreducible polynomials of $\mathbb{F}_p[x]$ of degree $f$.  Then $p$ {divides the index $i(K)$} if and only if $\mathcal{P}_f>\mathcal{N}_f$ for some positive integer $f$.
	\end{lem}

	{\it Proof of Theorem \ref{thm2}}.\\
	It is known that, in $\mathbb{F}_2[x]$, there are {two} monic irreducible polynomials of degree $1$ and just one monic irreducible polynomial of degree $2$, {namely} $x^2+x+1$. Since deg$(F)=6$, by Lemma \ref{index}, $2$  {divides the index $i(K)$}  if and only if there are at least {three} prime ideals of $\mathbb{Z}_K$ lying above $2$, with residue degree $1$ or at least {two} prime ideals of $\mathbb{Z}_K$ lying above $2$, with residue degree $2$.\\
	Since {$\Delta=-b^4(6^6b-5^5a^6)$ is the discriminant of $F(x)$, if $2$ divides $(\mathbb{Z}_K:\mathbb{Z}[\alpha])$, then $2$ divides $ab$}. Assume that $2$ divides $ab$, then
	\begin{enumerate}
		\item[(1)]If $2\mid a$ and $2\mid b$, then $F(x)\equiv x^6\md{2}$. Let $\phi=x$. Since $\nu_2(a)\geq1$, {then by assumption} $\nu_2(b)\leq 5$. {Thus} $N_{\phi}^+(F)=S$ has a single side
		joining $(0,\nu_2(b))\text{ and }(6,0)$. Three cases arise:
		\begin{enumerate}
			\item[(i)] If $\nu_2(b)\in\{1,5\}$, then the residual polynomial {of $F$} attached to $S$ is irreducible as it is of degree $1$. So, by Theorem $\ref{thm4}$, $2\mathbb{Z}_K=\mathfrak{p}^6$, with residue degree $1$.
			\item[(ii)] If $\nu_2(b)=3$, then the residual polynomial {of $F$} attached to $S$ is $R_{\lambda}(F)(y)=y^3+1=(y+1)(y^2+y+1)\in \mathbb{F}_{\phi}[y]$. So, by Theorem $\ref{thm4}$, $2\mathbb{Z}_K=\mathfrak{p}_{1}^2\mathfrak{p}_{2}^2$, with {$f_{1}=1\text{ and }f_{2}=2$}.
			\item[(iii)] If $\nu_2(b)\in\{2,4\}$, then the residual polynomial {of $F$} attached to $S$ is $R_{\lambda}(F)(y)=(y+1)^2\in \mathbb{F}_{\phi}[y]$, {which is not separable}. {Since $-\lambda=-1/3,-2/6$ is the slope of $S$ respectively, then $3$ divides the ramification index $e(\mathfrak{p})$ of any prime ideal $\mathfrak{p}$ of $\mathbb{Z}_{K}$ lying above $2$. Thus the possible cases are}:\\
			$2\mathbb{Z}_K=\mathfrak{p}^6$, with residue degree $1$.\\
			$2\mathbb{Z}_K=\mathfrak{p}^3$, with residue degree $2$.\\
			$2\mathbb{Z}_K=\mathfrak{p}_1^3\mathfrak{p}_2^3$, with residue degree $1$ each.
		\end{enumerate}
		\item[(2)]If $2\nmid a$ and $2\mid b$, then $F(x)\equiv x^5 (x+1)\md{2}$.  {Then $x+1$ provides a unique prime {ideal of $\mathbb{Z}_K$} lying above $2$, with residue degree $1$}. Let $\phi_1=x$.
		\begin{enumerate}
			\item[(i)] If $\nu_2(b)\not\equiv 0 \md{5}$, then $N_{\phi_1}^+(F)=S_{1}$ has a single side of degree $1$. Thus, by Theorem $\ref{thm4}$, $2\mathbb{Z}_K=\mathfrak{p}_{1}^5\mathfrak{p}_{2}$, with residue degree $1$ each.
			\item[(ii)] If $\nu_2(b)\equiv 0 \md{5}$, then $N_{\phi_1}^+(F)=S_{1}$ has a single side of degree $5$. Its attached residual polynomial {of $F$} is $R_{\lambda_1}(F)(y)=y^5+1=(y^4+y^3+y^2+y+1)(y+1)\in\mathbb{F}_{\phi_1}[y]$. Thus, by Theorem $\ref{thm4}$,  $2\mathbb{Z}_K=\mathfrak{p}_{11}\mathfrak{p}_{12}{\mathfrak{p}_{21}}$, 
			with {$f_{11}=4$ and $f_{12}=f_{21}=1$}.
		\end{enumerate}
		\item[(3)] If $2\mid a\text{ and }2\nmid b$, then $F(x)\equiv (x-1)^2(x^2+x+1)^2\md{2}$. {Let} $\phi_1=x-1$ and $\phi_2=x^2+x+1$, {then}
		\small
		$$\begin{array}{lll}
			F(x)&=&\phi_1^6+(a+6)\phi_1^5+(5a+15)\phi_1^4+(10a+20)\phi_1^3+(10a+15)\phi_1^2+(5a+6)\phi_1+a+b+1,\\
			&&\\
			&=&\phi_2^3+((a-3)x-2a)\phi_2^2+((a+2)x+3a-2)\phi_2-ax-a+b+1.
		\end{array}$$
		\normalsize
		{If $(a,b)\in \{(0,1),(2,3)\}\md{4}$, then by Theorem $\ref{thm1}$, $2\nmid (\mathbb{Z}_K:\mathbb{Z}[\alpha])$}. {Hence, $2\nmid i(K)$}.
		\begin{enumerate}
			\item[(i)] {For $(a,b)\equiv (0,3)\md{4}$},
			\begin{enumerate}
				\item[(a)] If $a\equiv4\md{8}$ and $b\equiv 3\md{8}$, then $N_{\phi_1}^+(F)=S_{11}+S_{12}$ has two sides joining the points $(0,u)$, $(1,1)$, and $(2,0)$, with $u\geq 3$. $N_{\phi_2}^+(F)=S_{21}$ has a single side joining $(0,2)\text{ and }(2,0)$ with the lattice point $(1,1)$ lying on it. Its attached residual polynomial {of $F$} is ${R_{\lambda_{21}}(F)(y)}=xy^2+(x+1)y+x$, which is irreducible over $\mathbb{F}_{\phi_2}$. By Theorem $\ref{thm4}$, $2\mathbb{Z}_K=\mathfrak{p}_{11}\mathfrak{p}_{12}\mathfrak{p}_{21}$, with {$f_{11}=f_{12}=1\text{ and }f_{21}=4$}. 
				\item[(b)]If $a\equiv4\md{8}$ and $b\equiv 7\md{8}$, then $N_{\phi_1}^+(F)=S_{1}$ has a single side joining $(0,2)\text{ and }(2,0)$ with the lattice point $(1,1)$ lying on it. Its attached residual polynomial {of $F$} is $R_{\lambda_1}(F)(y)=y^2+y+1$, which is irreducible over $\mathbb{F}_{\phi_1}$. $N_{\phi_2}^+(F)=S_{2}$ has a single side joining $(0,2)\text{ and }(2,0)$ with the lattice point $(1,1)$ lying on it. Its attached residual polynomial {of $F$} is $R_{\lambda_2}(F)(y)=xy^2+(x+1)y+x+1$, which is irreducible over $\mathbb{F}_{\phi_2}$. By Theorem \ref{thm4}, $2\mathbb{Z}_K=\mathfrak{p}_{1}\mathfrak{p}_{2}$, with {$f_{1}=2\text{ and }f_{2}=4$}.
				\item[(c)] If $a\equiv0\md{8}$ and $b\equiv 3\md{8}$, then $N_{\phi_1}^+(F)=S_{1}$ has a single side joining $(0,2)$ and $(2,0)$ with the lattice point $(1,1)$ lying on it. Its attached residual polynomial {of $F$} is $R_{\lambda_1}(F)(y)=y^2+y+1$, which is irreducible over $\mathbb{F}_{\phi_1}$. $N_{\phi_2}^+(F)=S_{2}$ has a single side joining $(0,2)$ and $(2,0)$ with the lattice point $(1,1)$ lying on it. Its attached residual polynomial {of $F$} is $R_{\lambda_2}(F)(y)=xy^2+(x+1)y+1=(y+1)(xy+1)\in\mathbb{F}_{\phi_2}[y]$. By Theorem \ref{thm4}, $2\mathbb{Z}_K={\mathfrak{p}_{11}}\mathfrak{p}_{21}\mathfrak{p}_{22}$, with residue degree $2$ each. {Hence, $2$ divides the index $i(K)$, and based on Engstrom's results \cite[page 234]{En}, $\nu_2(i(K))=2$}. 
				\item[(d)]If $a\equiv0\md{8}$ and $b\equiv 7\md{8}$, then $N_{\phi_1}^+(F)=S_{11}+S_{12}$ has two sides joining the points $(0,u),~(1,1)$, and $(2,0)$, with $u\geq 3$. $N_{\phi_2}^+(F)=S_{21}+S_{22}$ has two sides joining the points $(0,v),~(1,1)$, and $(2,0)$, with $v\geq 3$. Thus, by Theorem \ref{thm4}, $2\mathbb{Z}_K=\mathfrak{p}_{11}\mathfrak{p}_{12}\mathfrak{p}_{21}\mathfrak{p}_{22}$, with {$f_{11}=f_{12}=1\text{ and }f_{21}=f_{22}=2$}. Also in this case,  {$2$ divides the index $i(K)$, with $\nu_2(i(K))=2$}.
			\end{enumerate}
			\item[(ii)] For $(a,b)\equiv(2,1)\md{4}$,\\
			Since $\nu_2(a)=1$, $N_{\phi_2}^+(F)=S_2$ has a single side of height $1$. Thus, $\phi_2$ provides {a unique} prime {ideal of $\mathbb{Z}_K$} lying above $2$, with residue degree $2$ and ramification index $2$ either. {Let $u=\nu_2(5a+6)$ and $v=\nu_2(a+b+1)$}, {then we have the following cases}:
			\begin{enumerate}
				\item[(a)] If $v>2u$, then $N_{\phi_1}^+(F)=S_{11}+S_{12}$ has two sides joining the points $(0,v),~(1,u)$, and $(2,0)$. Thus, by Theorem \ref{thm4}, $2\mathbb{Z}_K=\mathfrak{p}_{11}\mathfrak{p}_{12}{\mathfrak{p}_{21}}^2$, with {${f}_{11}={f}_{12}=1$ and ${f_{21}}=2$}.
				\item[(b)] If $v=2u$; {$b\equiv -(a+1)+2^{2u}\md{2^{2u+1}}$}, then $N_{\phi_1}^+(F)={S_{1}}$ has a single side joining $(0,v)$ and $(2,0)$ with the lattice point $(1,u)$ lying on it. Its attached residual polynomial {of $F$} is $R_{\lambda_{1}}(F)(y)=y^2+y+1$, which is irreducible over $\mathbb{F}_{\phi_1}$. Thus, by Theorem $\ref{thm4}$, $2\mathbb{Z}_K={\mathfrak{p}_{1}}\mathfrak{p}_2^2$, with residue degree $2$ each.  {Hence, $2$ divides the index $i(K)$, and by Engstrom's results \cite[page 234]{En}, $\nu_2(i(K))=1$}.
				\item[(c)]If $v<2u$, then $N_{\phi_1}^+(F)=S_{1}$ has a single side joining $(0,v)$ and $(2,0)$.
				\begin{enumerate}
					\item[$(c_1)$] If $v$ is odd, then $S_{1}$ has degree $1$ and ramification index $2$, by Theorem \ref{thm4}, $2\mathbb{Z}_K={\mathfrak{p}_{1}}^2\mathfrak{p}_2^2$, with {${f_{1}}=1$ and $f_2=2$}.
					\item[$(c_2)$] If $v$ is even, then the residual polynomial {of $F$} attached to $S_{1}$ is $R_{\lambda_{1}}(F)(y)=(y+1)^2$, which is not separable over $\mathbb{F}_{\phi_1}$. Let $s$ be an integer such that $F(x)$ is $(x-s)$-regular with respect to $2$ and $s\equiv1\md{2}$, {then} $F(x)=\cdots+(15s^4+10as^3)(x-s)^2+(5as^4+6s^5)(x-s)+b+as^5+s^6$. {Since $\phi_2$ provides a unique prime ideal of $\mathbb{Z}_K$ lying above $2$, with residue degree $2$, we conclude that} $2$   {divides the index $i(K)$} if and only if $\nu_2(b+as^5+s^6) =2\nu_2(5as^4+6s^5)$. {In this case, $2\mathbb{Z}_K=\mathfrak{p}_1\mathfrak{p}_{2}^2$, with residue degree $2$ each. By \cite{En}, $\nu_2(i(K))=1$.}
				\end{enumerate}
			\end{enumerate}
		\end{enumerate}
	\end{enumerate}
	\begin{flushright}
		$\square$
	\end{flushright}
	{\it Proof of Theorem \ref{thm3}}.\\
	Since deg$(F)=6$, by Lemma \ref{index}, $3$ {divides the index $i(K)$} if and only if there are at least $4$ prime ideals of $\mathbb{Z}_K$ lying above $3$, with residue degree $1$.\\
	{Since $\Delta=-b^4(6^6b-5^5a^6)$, if $3$ divides the index} $(\mathbb{Z}_K:\mathbb{Z}[\alpha])$, then $3$ divides $ab$. Assume that $3$ divides $ab$, {then we have the following cases:}
	\begin{enumerate}
		\item[(1)] If $3\mid a\text{ and }3\mid b$, then $F(x)\equiv x^6\md{3}$. Let $\phi=x$. Since $\nu_3(a)\geq 1$, {by assumption} $\nu_3(b)\leq 5$. Thus $N_{\phi}^+(F)=S$ has a single side joining $(0,\nu_3(b))\text{ and }(6,0)$.
		\begin{enumerate}
			\item[(i)] If $\nu_3(b)\in\{1,5\}$, then $S$ is of degree $1$ and so, $3\mathbb{Z}_K=\mathfrak{p}^6$, with residue degree $1$.
			\item[(ii)] If $\nu_3(b)\in\{2,4\}$, then $S$ is of degree $2$. Its attached residual polynomial {of $F$} is $R_{\lambda}(F)(y)=y^2+b_3\in\mathbb{F}_{\phi}[y]$. If $b_3\equiv1\md{3}$, then $R_{\lambda}(F)(y)$ is irreducible over $\mathbb{F}_{\phi}$, by Theorem $\ref{thm4}$, $3\mathbb{Z}_K=\mathfrak{p}^3$, with residue degree $2$. If $b_3\equiv-1\md{3}$, then $R_{\lambda}(F)(y)=(y+1)(y-1)$ and so, $3\mathbb{Z}_K=\mathfrak{p}_{1}^3\mathfrak{p}_{2}^3$, with residue degree $1$ each.
			\item[(iii)] If $\nu_3(b)=3$, then $S$ has degree $3$, and its attached residual polynomial {of $F$} is $R_{\lambda}(F)(y)=(y+b_3)^3$, which is not separable over $\mathbb{F}_{\phi}$. {Since $-\lambda=-1/2$ is the slope of $S$, then $2$ divides the ramification index $e(\mathfrak{p})$ of any prime ideal $\mathfrak{p}$ of $\mathbb{Z}_K$ lying above $3$, and so} $\phi$ can provides at most $3$ prime ideals of $\mathbb{Z}_K$ lying above $3$, with residue degree $1$. {In this case, $3$ does not divide the index $i(K)$}.
		\end{enumerate}
		\item[(2)] If $3\nmid a \text{ and }3\mid b$, then $F(x)\equiv x^5(x+a)\md{3}$. Let $\phi_1=x$ and $\phi_2=x+a$, {then $\phi_2$ provides a unique prime {ideal of $\mathbb{Z}_K$} lying above $3$, with residue {degree} $1$}.  {On the other hand}, $N_{\phi_1}^+(F)=S_{1}$ has a single side joining $(0,\nu_3(b))$ and $(5,0)$. If $\nu_3(b)\not\equiv 0 \md{5}$, then  the degree of $S_{1}$ is $1$. Thus, by Theorem $\ref{thm4}$, $3\mathbb{Z}_K=\mathfrak{p}_{1}^5\mathfrak{p}_{2}$, with residue degree $1$ each. If $\nu_3(b)\equiv 0\md{5}$, then the residual polynomial {of $F$} attached to $S_1$ is $R_{\lambda_1}(F)(y)=a_3y^5+b_3=y^5+b_3a_3^{-1}$, which is separable over $\mathbb{F}_{\phi_1}$. If $b_3a_3^{-1}\equiv1\md{3}$, then $R_{\lambda_1}(F)(y)=(y+1)(y^4-y^3+y^2-y+1)$, and if $b_3a_3^{-1}\equiv-1\md{3}$, then $R_{\lambda_1}(F)(y)=(y-1)(y^4+y^3+y^2+y+1)$. It follows by Theorem $\ref{thm4}$ that $3\mathbb{Z}_K=\mathfrak{p}_{11}\mathfrak{p}_{12}\mathfrak{p}_{21}$, with {$f_{11}=f_{21}=1$ and $f_{12}=4$}.	
		\item[(3)] If $3\mid a\text{ and }3\nmid b$, then
		\begin{enumerate}
			\item[(i)] If $b\equiv 1 \md{3}$, then $F(x)\equiv(x^2+1)^3\md{3}$. Let $\phi=x^2+1$.
			{Since deg$(\phi)=2$, then $\phi$ cannot provide any prime ideal of $\mathbb{Z}_K$ lying above $3$, with residue degree $1$.}
			\item[(ii)] If $b\equiv -1\md{3}$, then $F(x)\equiv (x+1)^3 (x-1)^3\md{3}$. Let $\phi_1=x+1\text{ and }\phi_2=x-1$, {then}
			\small
			$$\begin{array}{lll}
				F(x)&=&\phi_1^6+(a-6)\phi_1^5+(-5a+15)\phi_1^4+(10a-20)\phi_1^3+(-10a+15)\phi_1^2+(5a-6)\phi_1-a+b+1,\\
				&&\\
				&=&\phi_2^6+(a+6)\phi_2^5+(5a+15)\phi_2^4+(10a+20)\phi_2^3+(10a+15)\phi_2^2+(5a+6)\phi_2+a+b+1.
			\end{array}$$
			\normalsize
			{If $(a,b)\in\{(0,2),(0,5),(3,-1),(6,-1)\}\md{9}$, then by Theorem $\ref{thm1}$, $3\nmid (\mathbb{Z}_K:\mathbb{Z}[\alpha])$}, {and so $3$ does not divide $i(K)$}.
			\item[(a)] {For $(a,b)\equiv(0,-1)\md{9}$, we have} $N_{\phi_i}^+(F)=S_{i1}+S_{i2}$ has two sides joining the points $(0,v_i),(1,1),\text{ and }(3,0)$, with $v_i\geq 3$ for every $i=1,2$. Thus, the degree of each side is $1$. It follows by Theorem $\ref{thm4}$ that $3\mathbb{Z}_K=\mathfrak{p}_{11}\mathfrak{p}_{12}^2\mathfrak{p}_{21}\mathfrak{p}_{22}^2$, with residue degree $1$ each. {Hence, $3$ divides the index $i(K)$, and by \cite{En}, $\nu_3(i(K))=1$}.
			
			\item[(b)] {For $(a,b)\equiv(3,5)\md{9}$, we have} $N_{\phi_1}^+(F)=S_{1}$ has a single side of height $1$. $N_{\phi_2}^+(F)=S_{21}+S_{22}$ has two sides joining the points $(0,v),~(1,1)\text{ and }(3,0)$, with $v\geq 2$. Thus, by Theorem $\ref{thm4}$, $3\mathbb{Z}_K={\mathfrak{p}_{11}}^3\mathfrak{p}_{21}\mathfrak{p}_{22}^2$, with residue degree $1$ each.
			
			\item[(c)] {For $(a,b)\equiv(6,5)\md{9}$}, we have $N_{\phi_2}^+(F)=S_{2}$ has a single side of height $1$ and $N_{\phi_1}^+(F)=S_{11}+S_{12}$ has two sides joining the points $(0,v),~(1,1)\text{ and }(3,0)$, with ${\tau}\geq 2$. Thus, by Theorem $\ref{thm4}$, $3\mathbb{Z}_K=\mathfrak{p}_{11}\mathfrak{p}_{12}^2\mathfrak{p}_{21}^3$, with residue degree $1$ each.
			
			\item[(d)] {For $(a,b)\equiv(3,2)\md{9}$, we have} $N_{\phi_2}^+(F)=S_2$ has a single side of height $1$ and ramification index $e(S_2)=3$. Let ${\tau}=\nu_3(-a+b+1)\geq 2\text{ and }{\mu}=\nu_3(5a-6)\geq 2$, then
			\begin{enumerate}
				\item[$(d_1)$] If ${\tau}=2$; $(a,b)\in\{(3,11),(3,20),(12,2),(12,20),(21,2),(21,11)\}\md{27}$, then $N_{\phi_1}^+(F)=S_{1}$ has a single side joining $(2,0)\text{ and }(3,0)$ of degree $1$ and ramification index $e(S_1)=3$.
				
				\item[$(d_2)$] If ${\tau}=3$ and ${\mu}=2$; $a\equiv 3,21\md{27}$ and $b\equiv a-1\pm 27\md{81}$, then $N_{\phi_1}^+(F)=S_{1}$ has a single side joining the points $(0,3)\text{ and }(3,0)$ with the lattice points $(1,2)$ and $(2,1)$ lying on it. Its attached residual polynomial {of $F$} is $R_{\lambda_1}(F)(y)=y^3+y^2+(5a-6)_3y+(-a+b+1)_3\in \mathbb{F}_{\phi_1}[y]$.
				
				\item[$(d_3)$] If ${\tau}=3$ and ${\mu}\geq3$; $a\equiv 12\md{27}$ and $b\equiv a-1\pm27\md{81}$ then $N_{\phi_1}^+(F)=S_{1}$ has a single side joining the points $(0,3)\text{ and }(3,0)$ with the lattice point $(2,1)$ lying on it. Its attached residual polynomial {of $F$} is $R_{\lambda_1}(F)(y)=y^3+y^2+(-a+b+1)_3\in \mathbb{F}_{\phi_1}[y]$.
				
				\item[$(d_4)$] If ${\tau}\geq4$ and ${\mu}=2$; $a\equiv 3,21\md{27}$ and $b\equiv a-1\md{81}$, then $N_{\phi_1}^+(F)=S_{11}+S_{12}$ has two sides joining the points $(0,{\tau}),~(1,2)\text{ and }(3,0)$, with the lattice point $(2,1)$ lying on $S_{12}$. The degree of $S_{11}$ is $1$ and the residual polynomial {of $F$} attached to $S_{12}$ is $R_{\lambda_{12}}(F)(y)=y^2+y+(5a-6)_3\in\mathbb{F}_{\phi_1}$.
				
				\item[$(d_5)$] If ${\tau}\geq4,~{\mu}\geq3$; $a\equiv 12\md{27}$ and $b\equiv a-1\md{81}$, then
				\begin{enumerate}
					\item[1.] If $2{\mu}={\tau}+1$, then $N_{\phi_1}^+(F)=S_{11}+S_{12}$ has a two sides joining the points $(0,{\tau}),(2,1)\text{ and }(3,0)$, with the lattice point $(1,u)$ lying on $S_{11}$. The degree of $S_{12}$ is $1$ and the residual polynomial {of $F$} attached to $S_{11}$ is $R_{\lambda_{11}}(F)(y)=y^2+(5a-6)_3y+(-a+b+1)_3\in \mathbb{F}_{\phi_1}[y]$.
					\item[2.] If $2{\mu}>{\tau}+1$, then $N_{\phi_1}^+(F)=S_{11}+S_{12}$ has a two sides joining the points $(0,{\tau}),(2,1)\text{ and }(3,0)$. The degree of $S_{12}$ is $1$. If ${\tau}$ is even, then the degree of $S_{11}$ is also $1$. Thus, $3\mathbb{Z}_K=\mathfrak{p}_{11}^2\mathfrak{p}_{12}{\mathfrak{p}_{21}}^3$, with residue degree $1$ each. If ${\tau}$ is odd, then the residual polynomial {of $F$} attached to $S_{11}$ is $R_{\lambda_{11}}(F)(y)=y^2+(-a+b+1)_3\in \mathbb{F}_{\phi_1}[y]$.
					\item[3.] If $2{\mu}<{\tau}+1$, then $N_{\phi_2}^+(F)=S_{11}+S_{12}+S_{13}$ has a three sides joining the points $(0,{\tau}),(1,{\mu}),(2,1)\text{ and }(3,0)$, with degree $1$ each. Thus, $3\mathbb{Z}_K=\mathfrak{p}_{11}\mathfrak{p}_{12}\mathfrak{p}_{13}{\mathfrak{p}_{21}}^3$, with residue degree $1$ each.
				\end{enumerate}
				{The following table (Table $1$) summarizes all cases treated in the case (d)}.
			\end{enumerate}
			\footnotesize
			\begin{table}[h]
				\centering
				\caption{}
				\begin{tabular}{|c|c|c|c|c|c|l}
					\cline{1-6}
					\multirow{3}{*}{${\tau}$}&\multirow{3}{*}{${\mu}$}&\multirow{2}{*}{$(-a+b+1)_3$}&\multirow{2}{*}{$(5a-6)_3$}&\multirow{3}{*}{$R_{\lambda}(F)(y)$}&\multirow{3}{*}{Q}&\\
					&&&&&&\\
					&&$\md{3}$&$\md{3}$&&&\\
					\cline{1-6}
					$2$&$\geq 2$&$\pm1$&&$y\pm1$&$2$&\\
					\cline{1-6}
					\multirow{6}{*}{$3$}&\multirow{4}{*}{$2$}& $1$&$1$& $(y+1)(y^2+1)$&$2$& \\
					\cline{3-6}
					&& $-1$&$1$& $y^3+y^2+y-1$&$1$& \\
					\cline{3-6}
					&& $1$&$-1$& $y^3+y^2-y+1$&$1$& \\
					\cline{3-6}
					&& $-1$&$-1$& $(y+1)^2(y-1)$& $\geq2$&$(1)$\\
					\cline{2-6}
					&\multirow{2}{*}{$\geq3$}&$1$&&$(y^2-y-1)(y-1)$&$2$ &\\
					\cline{3-6}
					&&$-1$&&$y^3+y^2-1$ &$1$&\\
					\cline{1-6}
					\multirow{15}{*}{$\geq 4$}&\multirow{2}{*}{$2$}&&$1$&$(y-1)^2$&$\geq2$&$(2)$\\
					\cline{3-6}
					&&&$-1$&$y^2+y-1$&$2$&\\
					\cline{2-6}
					&\multirow{2}{*}{${\mu}=\frac{{\tau}+1}{2}$}& $1$&$1$& $(y-1)^2$&$\geq 2$ &$(3)$\\
					\cline{3-6}
					&&$-1$&$1$&$y^2+y-1$&$2$&\\
					\cline{3-6}
					&\multirow{2}{*}{${\mu}\geq3$}&$1$&$-1$&$(y+1)^2$&$\geq 2$&$(4)$\\
					\cline{3-6}
					&&$-1$&$-1$&$y^2-y-1$&$2$&\\
					\cline{2-6}
					&{${\mu}>\frac{{\tau}+1}{2}$}&$1$&&$y^2+1$&$2$&\\
					\cline{3-6}
					&{${\mu}\geq 3$} &\multirow{2}{*}{$-1$}&&\multirow{2}{*}{$(y+1)(y-1)$}&\multirow{2}{*}{$4$}&\multirow{2}{*}{(5)}\\
					&$v$ is odd&&&&&\\
					\cline{2-6}
					&{${\mu}>\frac{{\tau}+1}{2}$}&&&&\multirow{3}{*}{$2$}&\\
					&{$u\geq3$}&&&&&\\
					&{$v$ is even}&&&&&\\
					\cline{2-6}
					&\multirow{2}{*}{${\mu}<\frac{{\tau}+1}{2}$}&\multirow{3}{*}{}&\multirow{3}{*}{}&&\multirow{3}{*}{$4$}&\multirow{3}{*}{(6)}\\
					&&&&&&\\
					&${\mu}\geq 3$&&&&\\
					\cline{1-6}
					\multicolumn{7}{l}{}\\
					\multicolumn{7}{l}{\textbf{Q: The number of prime ideals of $\mathbb{Z}_K$ lying above $3$, with residue degree $1$.}}
				\end{tabular}
			\end{table}
			\normalsize
			If the cases (1)-(4) (Table $1$) are not satisfied, then according to Table $1$, there are $4$ prime ideals of $\mathbb{Z}_K$ lying above $3$, with residue degree $1$ each if and only if one of the cases $(5)$ or $(6)$ arises. In this case, $3$ {divides the index $i(K)$} if and only if one of the cases $(5)$ or $(6)$ arises. i.e, $a\equiv 12\md{27}$, $b\equiv a-1\md{81}$, and $6\leq2{\mu}<{\tau}+1$ or $5\leq {\tau}+1<2{\mu},~{\tau}$ is odd, and $(-a+b+1)_3\equiv-1\md{3}$.\\
			{In the case $(5)$, we have $3\mathbb{Z}_K=\mathfrak{p}_{111}\mathfrak{p}_{112}\mathfrak{p}_{121}\mathfrak{p}_{211}^3$, with residue degree $1$ each. By \cite{En}, $\nu_3(i(K))=1$}.\\
			{In the case $(6)$, we have $3\mathbb{Z}_K=\mathfrak{p}_{11}\mathfrak{p}_{12}\mathfrak{p}_{13}\mathfrak{p}_{21}^3$, with residue degree $1$ each. By \cite{En}, $\nu_3(i(K))=1$}.\\
			
			Otherwise, if one of the cases (1)-(4) is satisfied, then we have: 
			\begin{enumerate}
				\item[$\bullet$] The case (1) (Table $1$) is satisfied if and only if $$(a,b)\in\{(21,74),(48,20),(75,47)\}\md{81}.$$
				Since $R_{\lambda_1}(F)(y)=(y+1)^2(y-1)$ is not separable over $\mathbb{F}_{\phi_1}$,						
				let $s$ be an integer such that $F(x)$ is $(x-s)$-regular with respect to $3$ and $s\equiv-1\md{3}$. Let $F(x)=\cdots+(20s^3+10as^2)(x-s)^3+(15s^4+10as^3)(x-s)^2+(5as^4+6s^5)(x-s)+b+as^5+s^6$. {Since $\phi_2$ provides a unique prime of $\mathbb{Z}_K$ lying above $3$, with residue degree $1$}, then $3$ {divides the index $i(K)$} if and only if $s$ satisfies one of the cases $(5)$ or $(6)$; 
				$5\leq\nu_3(b+as^5+s^6)+1<2\nu_3(5as^4+6s^5),~\nu_3(b+as^5+s^6)$ is odd, and $(b+as^5+s^6)_3\equiv-1\md{3}$ or
				$6\leq 2\nu_3(5as^4+6s^5)<\nu_3(b+as^5+s^6)+1$, {and we have $\nu_3(i(K))=1$}.
				
				\item[$\bullet$] The case $(2)$ is satisfied if and only if  $a\equiv3\md{27}$, and $b\equiv a-1\md{81}$. Since $R_{\lambda_1}(F)(y)=(y-1)^2$ is not separable over $\mathbb{F}_{\phi_1}$, let $s$ be an integer such that $F(x)$ is $(x-s)$-regular with respect to $3$ and $s\equiv-1\md{3}$, then $3$ is a {divides the index $i(K)$} if and only if $s$ satisfies one of the cases $(5)$ or $(6)$, {and we have $\nu_3(i(K))=1$}.
				
				\item[$\bullet$] The case $(3)$ is satisfied if and only if $a\equiv 12\md{27}$, $b\equiv a-1\md{81}$, $2{\mu=\tau}+1$, and $(-a+b+1)_3\equiv (5a-6)_3\equiv 1\md{3}$. Since $R_{\lambda_{11}}(F)(y)=(y-1)^2$ is not separable over $\mathbb{F}_{\phi_1}$, let $s$ be an integer such that $F(x)$ is $(x-s)$-regular with respect to $3$ and $s\equiv-1\md{3}$, then $3$ {divides the index $i(K)$} if and only if $s$ satisfies one of the cases $(5)$ and $(6)$, {and we have $\nu_3(i(K))=1$}. 
				
				\item[$\bullet$] The case $(4)$ is satisfied if and only if $a\equiv 12\md{27}$, $b\equiv a-1\md{81}$, $2{\mu=\tau}+1$, and $(-a+b+1)_3\equiv -(5a-6)_3\equiv 1\md{3}$. Since $R_{\lambda_{11}}(F)(y)=(y-1)^2$ is not separable over $\mathbb{F}_{\phi_1}$, let $s$ be an integer such that $F(x)$ is $(x-s)$-regular with respect to $3$ and $s\equiv-1\md{3}$, then $3$ {divides the index $i(K)$} if and only if $s$ satisfies one of the cases $(5)$ or $(6)$. {Also in this case, $\nu_3(i(K))=1$}.
			\end{enumerate} 
			
			\item[(e)] 
			For $(a,b)\equiv(6,2)\md{9}$. Then $N_{\phi_1}^+(F)=S_1$ has a single side of degree $1$ with ramification index $3$. Let ${v}=\nu_3(a+b+1)\geq 2$ and ${u}=\nu_3(5a+6)\geq2$, then
			\begin{enumerate}
				\item[$(e_1)$] If ${v}=2$; $(a,b)\in\{(6,2),(6,11),(15,2),(24,11),(24,20)\}\md{27}$, then $N_{\phi_2}^+(F)=S_2$ has a single side joining $(2,0)\text{ and }(3,0)$ of degree $1$ and ramification index $e(S_2)=3$.
				
				\item[$(e_2)$] If ${v}=3$ and ${u}=2$; $a\equiv 6,24\md{27}$ and $b\equiv -a-1\pm 27\md{81}$, then $N_{\phi_2}^+(F)=S_2$ has a single side joining the points $(0,3)\text{ and }(3,0)$ with the lattice points $(1,2)$ and $(2,1)$ lying on it. {Its attached residual polynomial {of $F$}} is $R_{\lambda_2}(F)(y)=-y^3+y^2+(5a+6)_3y+(a+b+1)_3\in \mathbb{F}_{\phi_2}[y]$.
				
				\item[$(e_3)$] If ${v}=3$ and ${u}\geq3$; $a\equiv 15\md{27}$ and $b\equiv -a-1 \pm 27\md{81}$, then $N_{\phi_2}^+(F)=S_2$ has a single side joining the points $(0,3)\text{ and }(3,0)$ with the lattice point $(2,1)$ lying on it. {Its attached residual polynomial {of $F$}} is $R_{\lambda_2}(F)(y)=-y^3+y^2+(a+b+1)_3\in \mathbb{F}_{\phi_2}[y]$.
				
				\item[$(e_4)$] If ${v}\geq4$ and ${u}=2$; $a\equiv 6,24\md{27}$ and $b\equiv-a-1\md{81}$, then $N_{\phi_2}^+(F)=S_{21}+S_{22}$ has two sides joining the points $(0,{v}),~(1,2)\text{ and }(3,0)$, with the lattice point $(2,1)$ lying on $S_{22}$. The degree of $S_{21}$ is $1$ and the residual polynomial {of $F$} attached to $S_{22}$ is $R_{\lambda_{22}}(F)(y)=-y^2+y+(5a+6)_3\in\mathbb{F}_{\phi_2}[y]$.
				
				\item[$(e_5)$]If ${v}\geq4,~{u}\geq3$; $a\equiv 15\md{27}$ and $b\equiv -a-1\md{81}$, then
				\begin{enumerate}
					\item[1.] If $2{u}={v}+1$, then $N_{\phi_2}^+(F)=S_{21}+S_{22}$ has a two sides joining the points $(0,{v}),(2,1)$, and $(3,0)$, with the lattice point $(1,{u})$ lying on $S_{21}$. The degree of $S_{22}$ is $1$ and the residual polynomial {of $F$} attached to $S_{21}$ is $R_{\lambda_{21}}(y)=y^2+(5a+6)_3y+(a+b+1)_3\in \mathbb{F}_2[y]$.
					
					\item[2.] If $2{u}>{v}+1$, then $N_{\phi_2}^+(F)=S_{21}+S_{22}$ has a two sides joining the points $(0,{v}),(2,1)$, and $(3,0)$. The degree of $S_{22}$ is $1$. If ${v}$ is even, then $S_{21}$ has degree $1$. Thus, $3\mathbb{Z}_K={\mathfrak{p}_{11}}^3\mathfrak{p}_{21}^2\mathfrak{p}_{22}$, with residue degree $1$ each. If ${v}$ is odd, then the residual polynomial {of $F$} attached to $S_{21}$ is $R_{\lambda_{21}}(y)=y^2+(a+b+1)_3\in \mathbb{F}_{\phi_2}[y]$.
					
					\item[3.] If $2{u}<{v}+1$, then $N_{\phi_2}^+(F)=S_{21}+S_{22}+S_{23}$ has a three sides joining the points $(0,{v}),(1,{u}),(2,1)\text{ and }(3,0)$, with degree $1$ each. Thus, $3\mathbb{Z}_K={\mathfrak{p}_{11}}^3\mathfrak{p}_{21}\mathfrak{p}_{22}\mathfrak{p}_{23}^3$, with residue degree $1$ each.
				\end{enumerate}
				{The following table (Table $2$) regroup all cases treated in the case (e).}
			\end{enumerate}
			\footnotesize
			\begin{table}[h]
				\centering
				\caption{}
				\begin{tabular}{|c|c|c|c|c|c|l}
					\cline{1-6}
					\multirow{3}{*}{${v}$}&\multirow{3}{*}{${u}$}&\multirow{2}{*}{$(a+b+1)_3$}&\multirow{2}{*}{$(5a+6)_3$}&\multirow{3}{*}{$R_{\lambda}(F)(y)$}&\multirow{3}{*}{Q}&\\
					&&&&&&\\
					&&$\md{3}$&$\md{3}$&&&\\
					\cline{1-6}
					$2$&$\geq 2$&$\pm1$&&$y\pm1$&$2$&\\
					\cline{1-6}
					\multirow{6}{*}{$3$}&\multirow{4}{*}{$2$}& $1$&$1$& $-y^3+y^2+y+1$&$1$& \\
					\cline{3-6}
					&& $-1$&$1$& $-(y-1)^2(y+1)$&$\geq 2$&$(1)$ \\
					\cline{3-6}
					&& $1$&$-1$& $-(y-1)(y^2+1)$&$2$& \\
					\cline{3-6}
					&& $-1$&$-1$& $-y^3+y^2-y-1$& $1$&\\
					\cline{2-6}
					&\multirow{2}{*}{$\geq3$}&$1$&&$-(y^2+y-1)(y+1)$&$2$ &\\
					\cline{3-6}
					&&$-1$&&$-y^3+y^2-1$ &$1$&\\
					\cline{1-6}
					\multirow{15}{*}{$\geq 4$}&\multirow{2}{*}{$2$}&&$1$&$-y^2+y+1$&$2$&\\
					\cline{3-6}
					&&&$-1$&$-(y+1)^2$&$\geq2$&$(2)$\\
					\cline{2-6}
					&\multirow{2}{*}{${u}=\frac{{v}+1}{2}$}& $1$&$1$& $(y-1)^2$&$\geq 2$ &$(3)$\\
					\cline{3-6}
					&&$-1$&$1$&$y^2+y-1$&$2$&\\
					\cline{3-6}
					&\multirow{2}{*}{${u}\geq3$}&$1$&$-1$&$(y+1)^2$&$\geq 2$&$(4)$\\
					\cline{3-6}
					&&$-1$&$-1$&$y^2-y-1$&$2$&\\
					\cline{2-6}
					&{${u}>\frac{{v}+1}{2}$}&$1$&&$y^2+1$&$2$&\\
					\cline{3-6}
					&{${u}\geq 3$} &\multirow{2}{*}{$-1$}&&\multirow{2}{*}{$(y+1)(y-1)$}&\multirow{2}{*}{$4$}&\multirow{2}{*}{(5)}\\
					&$\mu$ is odd&&&&&\\
					\cline{2-6}
					&${u}>\frac{{v}+1}{2}$&&&&\multirow{3}{*}{$2$}&\\
					&$\tau\geq 3$&&&&&\\
					&$\mu$ is even&&&&&\\
					\cline{2-6}
					&\multirow{2}{*}{${u}<\frac{{u}+1}{2}$}&\multirow{3}{*}{}&\multirow{3}{*}{}&&\multirow{3}{*}{$4$}&\multirow{3}{*}{(6)}\\
					&&&&&&\\
					&${u}\geq 3$&&&&\\
					\cline{1-6}
					\multicolumn{7}{l}{}\\
					\multicolumn{7}{l}{\textbf{Q: The number of prime ideals of $\mathbb{Z}_K$ lying above $3$, with residue degree $1$.}}
				\end{tabular}
			\end{table}
			\normalsize
			If the cases (1)-(4) (Table $2$) are not satisfied, then according to Table $2$, there are $4$ prime ideals of $\mathbb{Z}_K$ lying above $3$, with residue degree $1$ each if and only if one of the case $(5)$ or $(6)$ arises. In this case, $3$ {divides the index $i(K)$} if and only if one of the cases $(5)$ or $(6)$ arises. i.e, $a\equiv 15\md{27}$ and $b\equiv -a-1\md{81}$, and $6\leq2{u}<{v}+1$ or $5\leq{v}+1<2{u},~{v}$ is odd, and $(a+b+1)_3\equiv -1\md{3}$.\\
			{In the case $(5)$, we have $3\mathbb{Z}_K=\mathfrak{p}_{111}^3\mathfrak{p}_{211}\mathfrak{p}_{212}\mathfrak{p}_{221}$, with residue degree $1$ each. By \cite{En}, $\nu_3(i(K))=1$}.\\
			{In the case $(6)$, we have $3\mathbb{Z}_K=\mathfrak{p}_{11}^3\mathfrak{p}_{21}\mathfrak{p}_{22}\mathfrak{p}_{23}$, with residue degree $1$ each. Also in this case,  $\nu_3(i(K))=1$}.
			Otherwise, if one of the cases (1)-(4) is satisfied, then we have:
			\begin{enumerate}
				\item[$\bullet$] The case (1) (Table $2$) is satisfied if and only if $$(a,b)\in\{(6,47),(33,20),(60,74)\}\md{81}.$$ Since $R_{\lambda_2}(F)(y)=-(y-1)^2(y+1)$ is not separable over $\mathbb{F}_{\phi_2}$, let $s$ be an integer such that $F(x)$ is $(x-s)$-regular with respect to $3$ and $s\equiv1\md{3}$. Let $F(x)=\cdots+(20s^3+10as^2)(x-s)^3+(15s^4+10as^3)(x-s)^2+(5as^4+6s^5)(x-s)+b+as^5+s^6$. {Since $\phi_1$ provides a unique prime of $\mathbb{Z}_K$ lying above $3$, with residue degree $1$}, $3$ {divides the index $i(K)$} if and only if the regular $s$ satisfies one of the cases $(5)$ and $(6)$. i.e, if and only if 
				$5\leq\nu_3(b+as^5+s^6)+1<2\nu_3(5as^4+6s^5),~\nu_3(b+as^5+s^6)$ is odd, and $(b+as^5+s^6)_3\equiv-1\md{3}$ or
				$6\leq 2\nu_3(5as^4+6s^5)<\nu_3(b+as^5+s^6)+1$, {and we have $\nu_3(i(K))=1$}.
				
				\item[$\bullet$] The case $(2)$ is satisfied if and only if $a\equiv24\md{27}$, and $b\equiv -a-1\md{81}$. Since $R_{\lambda_2}(F)(y)=-(y+1)^2$ is not separable over $\mathbb{F}_{\phi_2}$, let $s$ be an integer such that $F(x)$ is $(x-s)$-regular with respect to $3$ and $s\equiv1\md{3}$, then $3$ {divides the index $i(K)$} if and only if $s$ satisfies one of the cases $(5)$ and $(6)$, {with $\nu_3(i(K))=1$}.
				
				\item[$\bullet$] The case $(3)$ is satisfied if and only if $a\equiv 15\md{27}$, $b\equiv -a-1\md{81}$, $2{u}={v}+1$, and $(a+b+1)_3\equiv(5a+6)_3\equiv 1\md{3}$. Since $R_{\lambda_{21}}(F)(y)=(y-1)^2$ is not separable over $\mathbb{F}_{\phi_2}$, let $s$ be an integer such that $F(x)$ is $(x-s)$-regular with respect to $3$ and $s\equiv1\md{3}$, then $3$ {divides index $i(K)$} if and only if $s$ satisfies one of the cases $(5)$ and $(6)$, {with $\nu_3(i(K))=1$}.
				
				\item[$\bullet$] The case $(4)$ is satisfied if and {only} if $a\equiv 15\md{27}$, $b\equiv -a-1\md{81}$, $2{u}={v}+1$, and $(a+b+1)_3\equiv-(5a+6)_3\equiv 1\md{3}$. Since $R_{\lambda_{21}}(y)=(y+1)^2$ is not separable over $\mathbb{F}_{\phi_2}$, let $s$ be an integer such that $F(x)$ is $(x-s)$-regular with respect to $3$ and $s\equiv1\md{3}$, then $3$ {divides the index $i(K)$} if and only if $s$ satisfies one of the cases $(5)$ and $(6)$, {with $\nu_3(i(K))=1$}.
			\end{enumerate}
		\end{enumerate}
	\end{enumerate}
	\begin{flushright}
		$\square$
	\end{flushright}

	{\it Proof of Proposition \ref{pro5}}.\\
	Since $\Delta=-b^4(6^6b-5^5a^6)$, if $5$	
	divides the index $(\mathbb{Z}_K:\mathbb{Z}[\alpha])$, then $5$ divides $b$. Assume that $5$ divides $b$, then
	\begin{enumerate}
		\item If $5\nmid a$, then $F(x)\equiv x^5(x+a)\md{5}$. {Let} $\phi_1=x$ and $\phi_2=x+a$, then $\phi_2$ provides a unique prime {ideal of $\mathbb{Z}_K$} lying above $5$, with residue degree $1$.
		\begin{enumerate}
			\item[(i)] If $\nu_5(b)\not\equiv 0\md{5}$, then $N_{\phi_1}^+(F)=S_1$ has a single side of degree $1$. Thus, by Theorem $\ref{thm4}$, $5\mathbb{Z}_K=\mathfrak{p}_{11}^5{\mathfrak{p}_{21}}$, with residue degree $1$ each. 
			
			\item[(ii)] If $\nu_5(b)\equiv 0\md{5}$, then
			$N_{\phi_1}^+(F)=S_1$ has a single side joining the points $(0,\nu_5(b))$ and $(5,0)$. Its attached residual polynomial {of $F$} is $R_{\lambda_1}(F)(y)=ay^5+b_5=(y+u)^5$, where $u=a^{-1}b_5\in \mathbb{F}_{\phi_1}^*$, which is not separable over $\mathbb{F}_{\phi_1}$. Let $\nu_5(b)=5k$, with $k$ is a positive integer.  Up to replace $\phi_1$ by $\phi_1+5^ku$, we take $\phi_1=x+5^ku$. Let $$F(x)=\phi_1^6+A_5\phi_1^5+A_4\phi_1^4+A_3\phi_1^3+A_2\phi_1^2+A_1\phi_1+A_0.$$
			With 
			$A_0=b-au^5(5^k)^5+u^6(5^k)^6$, 
			$A_1=5au^4(5^k)^4-6u^5(5^k)^5$, $A_2=15u^4(5^k)^4-10au^3(5^k)^3$, $A_3=10au^2(5^k)^2-20u^3(5^k)^3$, $A_4=15u^2(5^k)^2-5au(5^k)$, $A_5=a-6u5^k$.\\
			We have $\nu_5(A_0)\geq5k+1$, $\nu_5(A_1)\geq4k+1$, $\nu_5(A_2)=3k+1$, $\nu_5(A_3)=2k+1$, and $\nu_5(A_3)=k+1$.
			\begin{enumerate}
				\item[(a)] If $\nu_5(A_0)= 5k+1$, then $N_{\phi_1}^+(F)=S_{1}$ has a single side joining the points $(\nu_5(A_0))$ and $(5,0)$. Thus, $\phi_1$ provides {a unique} prime ideal of $\mathbb{Z}_K$ lying above $5$, with residue degree $1$.
				\item[(b)] If $\nu_5(A_0)\geq 5k+2$, then
				
				\item[$\bullet$] If $\nu_5(A_0)\leq 2\nu_5(A_1)-3k-1$,
				then $N_{\phi_1}^+(F)=S_{11}+S_{12}$ has two sides joining the points $(\nu_5(A_0)) $, $(2,3k+1)$, and $(5,0)$. Since the ramification index of $S_{12}$ is $e(S_{12})=3$, then $\phi_1$ can provides at most $3$ prime ideals of $\mathbb{Z}_K$ lying above $5$, with residue degree $1$ each. 
				
				\item[$\bullet$] If $\nu_5(A_0)>2\nu_5(A_1)-3k-1$, then $N_{\phi_1}^+(F)=S_{11}+S_{12}+S_{13}$ has three sides joining the points $(\nu_5(A_0)) $, $(1,\nu_5(A_1))$, $(2,3k+1)$, and $(5,0)$. Thus, $\phi_1$ provides $3$ prime ideals of $\mathbb{Z}_K$ lying above $5$, with residue degree $1$ each.
			\end{enumerate}
		\end{enumerate}
		In all these cases, there are at most $4$ prime ideals of $\mathbb{Z}_K$ lying above $5$, with residue degree $1$ each.
		\item If $5\mid a$, then $F(x)\equiv x^6\md{5}$. Since $\nu_5(a)\geq 1$, by assumption $v=\nu_5(b)\leq 5$. Thus, $N_{\phi_1}^+(F)=S$ has a single side joining $(0,v)$ and $(6,0)$, with ramification index $e(S)\in\{2,3,6\}$. So, at most there are $3$ prime ideals of $\mathbb{Z}_K$ lying above $5$, with residue degree $1$. 
	\end{enumerate}
	We conclude that $5$ {does not divide the index $i(K)$}.
	\begin{flushright}
		$\square$
	\end{flushright}
	\smallskip
	\begin{examples}
		Let $K$ be a  number field generated by a complex root of an irreducible trinomial  $F(x)\in \mathbb{Z}[x]$.
		\begin{enumerate}
			\item {For $F(x)=x^6+288x^5+154$, we have $(a,b)\equiv (0,1)\md{9}$, then by Theorem \ref{thm1} (iii), $\mathbb{Z}[\alpha]$ is not integrally closed. On the other hand, $(288,154)\equiv (0,10)\md{72}$, then by Corollary \ref{cor} (1), $i(K)=1$} .
			\item For $F(x)=x^6+18x^5+33$, we have $F(x) \equiv(x-1)^2(x^2+x+1)^2\md{2}$. By Theorem $\ref{thm2}$, $(x^2+x+1)$ provides a unique prime ideal of $\mathbb{Z}_K$ lying above $2$, with residue degree $2$ and  $F(x)$ is not $(x-1)$-regular. The regular integer in this case is $s=3$. Let $F(x)=\cdots+6075(x-3)^2+8748(x-3)+5136$. Thus, $N_{x-3}^+(F)=S$ has a single side joining $(0,4)$ and $(2,0)$, {with $R_\lambda(F)(y)=y^2+y+1$, which is irreducible over $\mathbb{F}_{x-3}$}. Then, $\mathbb{Z}_K=\mathfrak{p}_1\mathfrak{p}_2^2$, with residue degree $2$ each. Thus $2$ {divides $i(K)$, and by Theorem \ref{thm2} (3),  $\nu_2(i(K))=1$}. {On the other hand}, $(18,33)\equiv(0,6)\md{9}$, we conclude  by Theorem $\ref{thm3}$ that  $3$  does not divide  $i(K)$. Hence $i(K)=2$, and so $K$ is not monogenic.
			
			\item For $F(x)=x^6-42x^5-1258$, we have $F(x)\equiv(x+1)^3(x-1)^3\md{3}$. By Theorem $\ref{thm3}$, $(x-1)$ provides a unique prime ideal of $\mathbb{Z}_K$ lying above $3$, with residue degree $1$ and $F(x)$ is not $(x+1)$-regular. The regular integer in this case is $s=8\equiv-1\md{3}$, and we have $F(x)=\cdots-16640(x-8)^3-153600(x-8)^2-663552(x-8)-1115370$. Then $N_{x-8}^+(F)=S_{21}+S_{22}+S_{23}$ joining $(0,8)$, $(1,4)$, $(2,1)$, and $(3,0)$. Thus $3\mathbb{Z}_K={\mathfrak{p}_{11}}^3\mathfrak{p}_{21}\mathfrak{p}_{22}\mathfrak{p}_{23}$, with residue degree $1$ each. Thus $3$ {divides $i(K)$, and by Theorem \ref{thm3}, $\nu_3(i(K))=1$}.  {On the other hand, $(-42,-1258)\equiv (2,2)\md{4}$, then by Theorem $\ref{thm2}$, $2$ does not divide $i(K)$}. Hence  $i(K)=3$, and so $K$ is not monogenic.
			
			\item For $F(x)=x^6+144x^5+399$, {we have $(144,399)\equiv (0,39)\md{72}$, then by Corollary \ref{cor} (4), $i(K)=4$. Hence, $K$ is not monogenic}. 

			\item For $F(x)=x^6+54x^5+377$, we have $(54,377)\equiv (54,89)\md{288}$, then by Corollary \ref{cor} (5), $i(K)=6$. Hence, $K$ is not monogenic. 

			\item {For $F(x)=x^6+360x^5+35$, we have $(360,35)\equiv(0,35)\md{72}$, then by Corollary \ref{cor} (6), $i(K)=12$. Hence, $K$ is not monogenic}.

		\end{enumerate}
	\end{examples}
	\section*{Competing interests}
 There are  non-financial competing interests to report.

\end{document}